\documentclass[10pt]{amsart}

\usepackage{amsmath}  
\usepackage{amssymb}  
\usepackage{amsthm}   

\input xy
\xyoption{all}

\newtheorem{satz}{Theorem}[section]
\newtheorem{theorem}[satz]{Theorem}
\newtheorem{prop}[satz]{Proposition}
\newtheorem{lem}[satz]{Lemma}
\newtheorem{cor}[satz]{Corollary}
\theoremstyle{definition}
\newtheorem{defn}[satz]{Definition}
\newtheorem{rem}[satz]{Remark}

\newcommand{\cP}{\ensuremath{\mathcal{P}}}
\newcommand{\cM}{\ensuremath{\mathcal{M}}}
\newcommand{\cO}{\ensuremath{\mathcal{O}}}
\newcommand{\cF}{\ensuremath{\mathcal{F}}}
\newcommand{\cG}{\ensuremath{\mathcal{G}}}

\newcommand{\cZ}{\ensuremath{\mathcal{Z}}}

\newcommand{\ZZ}{\mathbb{Z}}
\newcommand{\NN}{\mathbb{N}}
\newcommand{\A}{\mathbb{A}}
\newcommand{\G}{\mathbb{G}}
\renewcommand{\P}{\mathbb{P}}

\newcommand{\Ab}{\textbf{Ab}}
\newcommand{\kalg}{\textbf{k-alg}}
\newcommand{\Spt}{\textbf{Spt}}
\newcommand{\Sm}{\textbf{Sm}}
\renewcommand{\dim}{\text{dim}}
\newcommand{\lci}{\text{l.c.i.}} 

\newcommand{\op}{\text{op}} 
\newcommand{\GL}{{\text{GL}}}

\DeclareMathOperator*{\colim}{colim}
\DeclareMathOperator*{\hocolim}{hocolim}
\DeclareMathOperator*{\hocofib}{hocofib}
\DeclareMathOperator*{\hofib}{hofib}
\DeclareMathOperator*{\Spec}{Spec}

\DeclareMathOperator*{\id}{id} 
\DeclareMathOperator*{\codim}{codim} 
\DeclareMathOperator*{\Hom}{Hom}

\author{Marc Levine, Christian Serp\'e}
\title{On a spectral sequence for equivariant K-theory}
\date{\today}   
\thanks{The first-named author gratefully acknowledges the support of the Humboldt Foundation through the Wolfgang Paul Program, and support of the NSF via grants DMS-0140445 and  DMS-0457195}

\begin{document}
\begin{abstract}
We apply the machinery developed by the first-named author to the $K$-theory of 
coherent $G$-sheaves on a finite type $G$-scheme $X$ over a field, where $G$ is a finite group. This leads to a definition of
$G$-equivariant higher Chow groups (different from the Chow groups of classifying spaces constructed by Totaro and generalized to arbitrary $X$ by  Edidin-Graham) and an Atiyah-Hirzebruch spectral sequence from the $G$-equivariant higher Chow groups to the higher $K$-theory of 
coherent $G$-sheaves on $X$. This spectral sequence generalizes the spectral sequence from motivic cohomology to $K$-theory constructed by Bloch-Lichtenbaum and Friedlander-Suslin.
\end{abstract}
\maketitle

\tableofcontents

\section{Introduction} There are many different approaches to a theory of equivariant algebraic $K$-theory; for simplicity, we restrict our discussion to the case of a finite group $G$ acting on a scheme $X$. One can use algebraic versions of the Borel construction, taking an sequence of finite dimensional approximations $U_n$ to $EG$ and then defining $K_*^G(X)$ as a limit of the $K_*((U_n\times X)/G)$. This approach gives a different answer depending on the underlying topology chosen; for the Zariski topology, one may use the standard simplicial model of $EG$ whereas for the \'etale topology it is necessary to use other models. This latter choice gives for example an equivariant $K_0$-group closely related to the equivariant Chow groups first defined by Totaro \cite{Totaro} in the case of $X$ a point, generalized to arbitrary $X$ by Edidin-Graham.\cite{EdidinGraham1}, and studied by Edidin-Graham \cite{EdidinGraham2}, Pandharipande \cite{Pandharipande} and many others.  One can also modify this approach by using the homotopy fixed-point spectrum $K(U_n\times X)^G$ instead of the $K$-theory of the quotient.

Another approach is to take the $K$-theory spectrum of the exact category of  locally free $G$-sheaves on $X$, or the $G$-theory analog using coherent $G$-sheaves. We will henceforth denote these two theories by $K_*(G,X)$ and $G_*(G,X)$. These theories differ in general from the ones described above, mainly due to the lack of \'etale descent for algebraic $K$-theory. 

It is possible to give an allied cycle theory by using the orbit category of $G$, as in topological Bredon (co)homology. This approach has been used by Joshua \cite{Joshua}, and has yielded for example Riemann-Roch type results for $G_*(G,X)$.

We will use a somewhat different approach. We feed the equivariant $G$-theory spectra $G(G,X)$ into the {\em homotopy coniveau} machine developed in \cite{TechLoc, LeHo}. What comes out is a theory of algebraic cycles with an interesting ``local coefficient" system, together with an associated theory of higher Chow groups and most importantly, a spectral sequence
\begin{equation}\label{eqn:SpecSeq}
E^1_{p,q}=CH_p(G,X,p+q)\Longrightarrow G_{p+q}(G,X).
\end{equation}
In the case $G=\{\id\}$, this reduces to the Bloch-Lichtenbaum/Friedlander-Suslin spectral sequence
\cite{BL, FS}.

The local coefficient system is easy to describe: Let $Y$ be a finite type $k$-scheme with a $G$-action, $Z\subset Y$ an irreducible closed subset with generic point $z$. Let $G_z\subset G$ be the isotropy group of $z$. Then $G_z$ acts on the residue field $\kappa(z)$, and we define the ``coefficient group" of $Z$ to be $K_0(G_z, \kappa(z))$. The dimension $p$ equivariant cycle group of $Y$ is then
\[
z_p(G,Y):=\oplus_{[z]\in (Y_{(p)})/G}K_0(G_z,\kappa(z)).
\]
To explain: $Y_{(p)}$ is the set of points $z\in Y$ with closure $\bar{z}$ having dimension $p$ over $k$. $G$ acts on $Y_{(p)}$ and $z$ is a representative of the orbit $[z]\in (Y_{(p)})/G$.

Now if $X$ is a finite type $k$-scheme with $G$-action, taking the dimension $p+r$ points of $X\times\Delta^r$ in ``good position" with respect to the faces of $\Delta^r$, as in Bloch's construction of his cycle complexes \cite{AlgCyc}, and using our coefficient system $[z]\mapsto K_0(G_z,\kappa(z))$, one defines the Bredon-type equivariant cycle complex $z_p(G,X,*)$ and equivariant higher Chow groups $CH_p(G,X,*)$. The arguments of \cite{TechLoc, LeHo} go through without much change to yield the spectral sequence \eqref{eqn:SpecSeq}.

As alluded to above, one can view our construction as giving a Bredon-type motivic Borel-Moore homology theory for the quotient stack $[X/G]$. It seems reasonable that one could extend the construction of the cycle theory and the spectral sequence to more general stacks, either by removing the restriction that $G$ is finite (e.g. allowing $G$ to be a linear algebraic group) or considering stacks other than quotient stacks. It would also be interesting to see if Joshua's  construction, applied to the Friedlander-Suslin-Voevodsky cycle complexes, yield the same motivic Borel-Moore homology groups as ours. 

Part of the development of this paper was made possible by the hospitality of the University of Duisberg-Essen, the University of M\"unster and  the Wolgang Paul program of the Alexander von Humboldt Foundation. We would like to express our sincere gratitude for this support.

\vspace{2cm}

\section{Equivariant $G$ and $K$-theory}
Let $G$ be a finite group and $X$ be a (left) $G$-scheme. A coherent $G$-module $\cF$ on $X$
is a coherent $\cO_X$-module together with an action of $G$ on this
module which is compatible with the action of $G$ on $X$. 
In other words, for each $g\in G$ we have a morphism
$\phi_g:g^*\cF\rightarrow \cF$ and for $g,h\in G$ we have
$\phi_{g\cdot h}=\phi_h\circ h^* \phi_g$.
We denote by $\cM_{G,X}$ 
the abelian category of coherent $G$-modules and by
$\cP_{G,X}$ the exact subcategory of $\cM_{G,X}$ of those
G-modules which are locally free as $\cO_X$-modules. We further denote by
 $G(G,X)$ and $K(G,X)$ the K-theory spectrum of $\cM_{G,X}$ resp. $\cP_{G,X}$.
Then $G(G,X)$ is contravariant with respect to flat $G$-morphisms in $X$
and covariant with respect to projective $G$-morphisms in $X$.
$K(G,X)$ is contravariant with respect to any $G$-morphism.

If $X=\Spec R$ for a commutative ring $R$ the (left) action of $G$
on $X$ induces a (right) action of $G$ on the ring $R$. We set
$R^{tw}[G]:=\oplus_{g\in G} Re_g$. By 
$$(r_g\cdot e_g)(r_h\cdot e_h):=r_g\cdot (r_h\cdot g^{-1}) e_{(g\cdot h)}$$
for all $g,h\in G$ and $r_g,r_h\in R$ we define a ring structure
on $R^{tw}[G]$ and we call $R^{tw}[G]$ the twisted group ring 
of $R$.

\begin{lem}\label{twiequ}
Let $X=\Spec R$ be a noetherian affine $G$-scheme.
Then the  category of finitely generated $R^{tw}[G]$-modules is 
equivalent to $\cM_{G,X}$.
\end{lem}
\begin{proof}
For a  $R$-module $M$ and $g\in G$ let $\;_gM$ be the $R$-module
which is $M$ as abelian group and where the $R$-module structure
is defined by $r\cdot_{\;_gM}m:=(r\cdot g^{-1})\cdot m$ for all $r\in R$ and
$m\in M$. A $G$-module on $\Spec R$ is the same as an $R$-module
together with $R$-module homomorphisms
$\phi_g:\;_gM\rightarrow M$ for all $g\in G$ with 
$\phi_{(gh)}=\phi_h\circ \phi_g$ for all $g,h\in G$.
Then $M$ becomes a left $R^{tw}[G]$-module if we set
$e_g\cdot m := \phi_{g^{-1}}(m).$
Conversely if $M$ is a $R^{tw}[G]$-module we define
$\phi_g:\;\;_gM\rightarrow M$ by $\phi_g(m):=e_{g^{-1}}\cdot m$
for all $g\in G$ and get in this way a $G$-module
on $\Spec R$. Obviously the two functors are inverse to other. 
Because $G$ is finite and $R$ noetherian the lemma follows.
\end{proof}

\begin{lem}\label{lemmaequik}
Let $R$ be commutative noetherian ring with $\frac{1}{\sharp G}\in R$. Then the category
of finitely generated projective $R^{tw}[G]$-modules is equivalent
to $\cP_{G,\Spec R}$.
\end{lem}
\begin{proof}
If $\sharp G$ is invertible in $R$ then each $R^{tw}[G]$-module which
is projective as $R$-module is also projective as
$R^{tw}[G]$-module. So the lemma follows from the lemma above.
\end{proof}

\vspace{0.5cm}

Now let $k$ be a fixed commutative ring with $\frac{1}{\sharp G}\in k$. We denote by $\kalg$ the category
of commutative $k$-algebras. Let $F:\kalg\rightarrow \Ab$ be an additive functor.
Following Vorst \cite{Vorst}, we define
$$
\begin{array}{cccc}
  NF: & \kalg & \rightarrow & \Ab \\
     & R   & \mapsto     & ker(F(R[T])\rightarrow F(R)),
\end{array}
$$
where the last map is induced by $R[T]\rightarrow R, T\mapsto 0.$
For $q>1$ we define inductively  $N^qF:=N(N^{q-1}F)$.

For $f\in R$ the morphism $R[X]\rightarrow R[X], X\mapsto f\cdot X$
induces a group endomorphism $NF(R)\rightarrow NF(R)$. So $NF(R)$ becomes
a $\ZZ[T]$ module. We denote by $NF(R)_{[f]}$ the $\ZZ[T,T^{-1}]$ module
$\ZZ[T,T^{-1}]\otimes_{\ZZ[T]} NF(R)$. With these notations
Vorst proves the following theorem.
\begin{theorem}\label{theorem23}
Let $R\in \kalg$ and let $r_1,\dots, r_n$ be elements of $R$ which generate
the unit ideal. Suppose further that the map
$$NF(R[T]_{r_{i_0},\dots,\hat{r_{i_j}},\dots,r_{i_p}})_{[r_{i_j}]}\rightarrow
                 NF(R[T]_{r_{i_0},\dots,r_{i_p}})$$
is an isomorphism, for each set of indexes $1\leq i_0<\dots <i_p\leq n$. Then
the canonical morphism
$$\epsilon: NF(R)\rightarrow \oplus_{j=1}^n NF(R_{r_j})$$
is injective.
\end{theorem}
\begin{proof}
Compare \cite[Theorem 1.2]{Vorst} or \cite[Lemma 1.1]{LevineRevisited}.
\end{proof}

As in \cite{TT}, one may form the {\em Bass delooping} $KB(G,X)$ of the -1-connected spectrum $K(G,X)$. The functor $X\mapsto KB(G,X)$ satisfies
\begin{enumerate}
\item There is a canonical map $K(G,X)\to KB(X,G)$, identifying $K(G,X)$ with is the -1-connected cover of $KB(G,X)$.
\item For every $X$, there is the natural exact sequence
\begin{multline*}
0\to KB_p(G,X)\to KB_p(G,X\times\A^1)\oplus KB_p(G,X\times\A^1)\\
\to KB_p(G,X\times\G_m)\to
KB_{p-1}(G,X)\to 0
\end{multline*}
called the {\em fundamental exact sequence}.
\item If $X$ is regular, then $K(G,X)\to KB(G,X)$ is a weak equivalence.
\end{enumerate}
From now on, we will drop the notation $KB(G,X)$ and write $K(G,X)$ for the (possibly) non-connected version.

Now let $X=\Spec C$ be an affine $G$-scheme over $k$. We consider for each $p\in\NN$ the
functor
$$
\begin{array}{cccc}
  K_p(G,X\otimes_k -): & \kalg & \rightarrow & \Ab \\
     & R   & \mapsto     & K_p(G,X\otimes_k R),
\end{array}
$$
where $G$ acts trivially on $R$. We say that $X$ 
is $K_p(G,-)$-regular if $N^qK_p(G,X)=0$
for all $q>0$.

\begin{lem}\label{lemmaNp}
Let $R\in \kalg$ and $f\in R$. Suppose that there is an $g\in R$ such that $fg=0$ and
$f+g$ is a  non-zero divisor. Then the natural map
$$N^qK_p(G,X\otimes_k R)_{[f]}\rightarrow N^qK_p(G,X\otimes_k R_f)$$
is an isomorphism.
\end{lem}
\begin{proof}
By Lemma \ref{lemmaequik} we have the identification
$$K_p(G,X\otimes_k R)=K_p( (C\otimes_k R)^{tw}[G]).$$
Furthermore the element $f\in C\otimes_k R$ lies in the
center of the ring. The proof of  Lemma 1.4 in \cite{Vorst}
then goes through word for word to prove the lemma.
\end{proof}

\begin{theorem}\label{theoremreg}
Let $R\in \kalg$ be reduced and $f_1,\dots,f_n \in R$ such that
$(f_1,\dots,f_n)=R$ and such that $X\otimes_k R_{f_i}$
is $K_p(G,-)$-regular for $i=1,\dots, n$. Then
$X\otimes_k R$ is also $K_p(G,-)$-regular.
\end{theorem}
\begin{proof}
This follow from theorem \ref{theorem23} and lemma \ref{lemmaNp}.
\end{proof}

\begin{cor}\label{correg}
Let $X=\Spec C$ be a reduced affine $G$-scheme of finite type over a field and let $U_i\subset X$, $i=1,\dots,n$ 
an affine $G$-invariant (i.e. $GU_i\subset U_i$) open covering of $X$. Then $X$ is $K_p(G,-)$-regular if and only 
if $U_i$ is $K_p(G,-)$-regular for $i=1,\dots, n$.
\end{cor}
\begin{proof}
Let $X/G:=\Spec C^G$ be the quotient $X$ by $G$ and $\pi:X\rightarrow X/G$ the projection. Then $\pi(U_i)$ are open
subsets of $X/G$ and since the $U_i$'s are $G$-invariant we have $\pi^{-1}\pi(U_i)=U_i$ (see for example \cite{SGA1}[V.1.1]). 
For each $i=1,\dots, n$ we write $\pi(U_i)$ as a union of the principle open subschemes
\[
\pi(U_i)=\cup_j (X/G)_{f_{ij}}
\]
for suitable $f_{i1}, \dots, f_{im_i}\in C^G$. Noting that 
$(X/G)_{f_{ij}}=\Spec(C^G_{f_{ij}})$ and letting  $U_{ij}:=X\otimes_{C^G}(C^G_{f_{ij}})$, we have $U_i=\bigcup_j U_{ij}$.
We apply lemma \ref{lemmaNp} and theorem \ref{theoremreg} with $k=C^G$ to get the result.

\end{proof}
\vspace{0.5cm}

In the last section of this paper we need to consider the equivariant $K$-theory of some singular schemes. This can be partially understood by comparing with multi-relative $K$-theory and $KH$-theory, which we now recall.

For an $G$-invariant subscheme $Y$ of a scheme $X$ we define the relative
equivariant K-theory spectrum as
$$K(G,X;Y):=\hofib(K(G,X)\rightarrow K(G,Y)).$$
More generally,  for a family $\{Y_1,\dots ,Y_n\}$
of $G$-invariant subschemes we inductively define
\begin{multline*}
K(G,X;Y_1,\dots ,Y_n):=\\
\hofib[K(G,X;Y_1,\dots,Y_{n-1})\rightarrow 
K(G,Y_n;Y_1\cap Y_n,\dots,Y_{n-1}\cap Y_n)].
\end{multline*}
It can easily be seen that this definition is independent
of the order of $\{Y_1,\dots,Y_n\}$.

\vspace{0.5cm}
For an affine $G$-scheme $X=\Spec R$ we denote by
$$KH(G,X):=KH(R^{tw}[G])$$
where on the right hand side $KH$ denotes the homotopy
K-theory of Weibel for (not necessarily commutative!) rings (compare \cite{WeibelHomo}).

If $\{Y_1,\dots ,Y_n\}$ is a family of $G$-invariant affine
subschemes of the affine $G$-scheme $X$ we define
as above $KH(G,X;Y_1,\dots, Y_n)$. 

There is a natural transformation of functors from rings to spectra
\[
K(-)\to KH(-);
\]
this induces a natural transformation of functors from rings with a $G$-action to spectra
\[
K(G,-)\to KH(G,-).
\]
Since $K(G,T)\to KH(G,T)$ is a weak equivalence for $T$ regular, the map
\[
K(G,X;Y_1,\dots, Y_n)\to KH(G,X;Y_1,\dots, Y_n)
\]
is a weak equivalence if $X$ and all the intersections $Y_{i_1}\cap\ldots\cap Y_{i_s}$ are regular.

Also, if $Y=\cup_{i=1}^nY_i$, then it is easy to see that
\[
K(G,X;Y)=K(G,X;Y,Y)=\ldots=K(G,X;Y,\ldots,Y).
\]
and similarly for $KH(G,-)$. The inclusions $Y_i\to Y$ thus induce the maps
\[
\alpha:K(G,X;Y)=K(G,X;Y,\ldots,Y)\to K(G,X;Y_1,\ldots, Y_n)
\]
and 
\[
\beta:KH(G,X;Y)=KH(G,X;Y,\ldots,Y)\to KH(G,X;Y_1,\ldots, Y_n).
\]
By from \cite[th. 1.3,Cor. 2.2]{WeibelHomo},  $\beta$ is a weak equivalence. 

For the reader's convenience, we include a proof of the following elementary fact:

\begin{lem}\label{lem:K1MV} Let $S$ be a    ring  (not necessarily commutative). Let $I_1, I_2\subset S$ be two-sided ideals, and let $S_i:=S/I_jS$, $j=1,2$ and $S_{12}:=S/(I_1+I_2)S$. Let $\pi_i:S\to S_i$,
$\pi_{12,i}:S_i\to S_{12}$ be the quotient maps and let 
\[
S_1\times^{S_{12}}S_2=\{(s_1, s_2)\in S_1\times S_2\ |\ \pi_{12,1}(s_1)=\pi_{12,2}(s_2)\}. 
\]
Suppose that   the map
\begin{align*}
&\pi:S\to S_1\times_{S_{12}}S_2\\
&\pi(s):=(\pi_1(s), \pi_2(s))
\end{align*}
is an isomorphism. Suppose in addition that the surjection $S_2\to S_{12}$ is split by a ring homomorphism $\sigma_2:S_{12}\to S_2$. Then the sequence
\[
0\to K_1(S)\xrightarrow{(\pi_{1*},\pi_{2*})} K_1(S_1)\times K_1(S_2)\xrightarrow{\pi_{12,1*}-\pi_{12,2*}}
K_1(S_{12})\to0
\]
is  exact.
\end{lem}

\begin{proof}  We consider the $K_1$-sequence. We first show that $(\pi_{1*},\pi_{2*})$ is injective. If $\alpha\in \GL_N(S)$ goes to zero in $K_1(S_1)\times K_1(S_2)$ then we can write
\begin{align*}
&\pi_1(\alpha)=\prod_{r=1}^n e_{i_r,j_r}^{\lambda_r}\in \GL_M(S_1),\\
&\pi_2(\alpha)=\prod_{r=1}^n e_{i_r,j_r}^{\mu_r}\in \GL_M(S_2)
\end{align*}
 for some $M\ge N$ and for some $\lambda_r\in S_1$, $\mu_r\in S_2$ (we may take some $\lambda_r$ or $\mu_s$ zero if need be). Letting $\bar{\lambda}_r=\pi_{12,1}(\lambda_r)$    and 
  $\bar{\mu}_r=\pi_{12,1}(\mu_r)$, we have the relation
\[
\prod_{r=n}^1 e_{i_r,j_r}^{-\bar{\mu}_r}\prod_{r =1}^ne_{i_r,j_r}^{\bar{\lambda}_r}=1
\]
in $\GL_M(S_{12})$. Thus the element of the Steinberg group $St_M(S_{12})$ given by 
\[
\prod_{r=n}^1 x_{i_r,j_r}^{-\bar{\mu}_r}\prod_{r =1}^nx_{i_r,j_r}^{\bar{\lambda}_r}
\]
 defines an element $x\in K_2(S_{12})$. Lift $x$ to the element 
 \[
 x_2:=\sigma_*(x)=\prod_{r=n}^1 x_{i_r,j_r}^{-\tilde{\mu}_r}\prod_{r =1}^nx_{i_r,j_r}^{\tilde{\lambda}_r}
\]
in $K_2(S_2)$; here $\tilde{\lambda}_r=\sigma(\bar{\lambda}_r)$ and similarly for $\tilde{\mu}_r$.

Lift each $\bar{\mu}_r$ to an element $\rho_r\in S_1$. Then
\begin{align*}
\pi_1(\alpha)=\prod_{r=1}^n e_{i_r,j_r}^{\rho_r}\prod_{r=n}^1 e_{i_r,j_r}^{-\rho_r}\prod_{r=1}^n e_{i_r,j_r}^{\lambda_r}\\
\pi_2(\alpha)=\prod_{r=1}^n e_{i_r,j_r}^{\mu_r}\prod_{r=n}^1 e_{i_r,j_r}^{-\tilde{\mu}_r}\prod_{r =1}^ne_{i_r,j_r}^{\tilde{\lambda}_r}
\end{align*}
Clearly the pairs $(\rho_r,\mu_r)$, $(-\rho_r,-\tilde{\mu}_r)$, $(\lambda_r\tilde{\lambda}_r)$ define elements $a_r$, $b_r$ and $c_r$ in $S_1\times^{S_{12}}S_2=S$ and we have
\[
\alpha=\prod_{r=1}^n e_{i_r,j_r}^{a_r}\prod_{r=n}^1 e_{i_r,j_r}^{b_r}\prod_{r=1}^n e_{i_r,j_r}^{c_r},
\]
as desired.

The exactness at $K_1(S_1)\times K_1(S_2)$ is easier: if $\pi_{12,1*}(\alpha_1)=
\pi_{12,2*}(\alpha_2)$ in $K_1(S_{12})$ for $\alpha_i\in \GL_{N_i}(S_i)$, $i=1,2$, then there is an element $e\in E_N(S_{12})$ with
\[
\pi_{12,1*}(\alpha_1)=\pi_{12,2*}(\alpha_2)e
\]
in $\GL_{N}(S_{12})$ for some $N\ge N_1, N_2$. We can lift $e$ to an $e_2\in E_N(S_2)$; replacing $\alpha_2$ with $\alpha_2e_2$, we may assume that  $N_1=N_2=N$ and
\[
\pi_{12,1*}(\alpha_1)=\pi_{12,2*}(\alpha_2)
\]
in $\GL_{N}(S_{12})$. Since $S=S_1\times^{S_{12}}S_2$, there is a unique $\alpha\in \GL_N(S)$ with $\pi_i(\alpha)=\alpha_i$, $i=1,2$.

The exactness at $K_1(S_{12})$ follows by using the splitting $\sigma_*$ to $\pi_{12,2*}$.
\end{proof}

Let 
\[
\Delta^n:=\Spec k[t_0,\ldots, t_n]/\sum_it_i-1.
\]
Let $\partial_i\Delta^n$ be the closed subscheme of $\Delta^n$ defined by $t_i=0$ and $\partial\Delta^n=\cup_{i=0}^n\partial_i\Delta^n$, . For a $k$-scheme with $G$-action $X$, we have the product schemes $X\times\Delta^n$, $X\times\partial\Delta^n$, with $G$ acting by the identity on $\Delta^n$, $\partial\Delta^n$, and by the given action on $X$. 

\begin{lem}\label{lem:Regularity}
Let $k$ be a field, $G$ a finite group with $\frac{1}{\sharp G}\in k$,
$X$ a regular affine $G$-scheme over $k$, and let $U$ be a $G$-stable affine open subscheme of $X\times\partial\Delta^n$. Then  $U$ is $K_p(G,-)$-regular for all $p\le 1$.
\end{lem}

\begin{proof} It suffices to show that $U$ is $K_1$-regular: we have for each $p$ and  $q$ the fundamental exact sequence
\begin{multline*}
0\to  N^qK_p(G,U)\\
 \to N^qK_p(G,U\times_k\A^1)\oplus N^qK_p(G,U\times_k\A^1)\\\to
N^qK_p(G,U\times_k\G_m)\to N^qK_{p-1}(G,U)\to0
\end{multline*}
where $G$ acts trivially on $\A^1$, $\G_m$. Clearly, if $U$ is $K_p(G,-)$ regular, so is $U\times\A^1$; by Lemma~\ref{lemmaNp}, this implies $U\times\G_m$ is also $K_p(G,-)$ regular. The exact sequence shows that $U$ is $K_{p-1}(G,-)$ regular.

Let $\partial_{\le m}\A^n$ be the closed subscheme of $\Spec k[x_1,\ldots, x_n]$ defined by $\prod_{i=1}^mx_i=0$. Note that each point of $U$ has an open neighborhood isomorphic (as a $G$-scheme) to an open neighborhood of $X\times\partial_{\le m}\A^n$ for some $m$. 
By Corollary \ref{correg}, it suffices to show that  $X\times\partial_{\le m}\A^n$ is $K_1(G,-)$ regular. We prove this by induction on $n$ and $m$, the case of arbitrary $n$ and $m=1$ following from the smoothness of $\partial_{\le 1}\A^n=\A^{n-1}$ over $k$.

Write $R_{n,m}=k[x_1,\ldots, x_n]/\prod_{i=1}^mx_i$. Note that 
\[
\partial_{\le m}\A^n=\partial_{\le m-1}\A^n\cup_{\partial_{\le m-1}\A^{n-1}}\A^{n-1}.
\]
and the inclusion $\partial_{\le m-1}\A^{n-1}\to \partial_{\le m-1}\A^n$ is split by the projection
\[
(x_1,\ldots, x_n)\mapsto (x_1,\ldots,x_{m-1},x_{m+1},\ldots, x_n),
\]
Thus, for each affine $G$-scheme $Y=\Spec A$ over $k$, we have  the ring
\[
S:=A^{tw}[G]\otimes_kR_{n,m}
\]
with two-sided ideals $I_1=(\prod_{i=1}^{m-1} x_i)$, $I_2=(x_m)$ which satisfy the hypotheses of Lemma~\ref{lem:K1MV}. Note also that $S/I_1=A^{tw}[G]\otimes_kR_{n,m-1}$ and $S/I_2=A^{tw}[G]\otimes_kR_{n,1}$, so the sequence
\begin{multline*}
0\to K_1(G,Y\times\partial_{\le m}\A^n)\\
\to K_1(G,Y\times\partial_{\le m-1}\A^n)\oplus K_1(G,Y\times\A^{n-1})\\\to
K_1(G,Y\times\partial_{\le m-1}\A^{n-1})\to0
\end{multline*}
is exact. Taking $Y=X\times\A^q$ and using our induction hypothesis shows that $X\times\partial_{\le m}\A^n$ is $K_1$-regular.
\end{proof}

\begin{prop}\label{homfibseq}
Let $k$ be a field, $G$  a finite group with $\frac{1}{\sharp G}\in k$,
$X$   a regular affine $k$-scheme with a $G$-action, and  $U$   a $G$-stable affine open subscheme of $X\times\Delta^n$. Let $\partial U:=X\times\partial\Delta^n\cap U$, $\partial_iU:= U\cap X\times\partial_i\Delta^n$. Then the map
\[
\alpha:K_0(G,U;\partial U)\to K_0(G,U;\partial U_0,\ldots, \partial U_n)
\]
is an isomorphism.
\end{prop}
\begin{proof} Both $U$ and all the intersections $\partial U_{i_1}\cap\ldots\cap\partial U_{i_s}$ are regular, so 
\[
K(G,U;\partial U_0,\ldots, \partial U_n)\to
KH(G,U;\partial U_0,\ldots, \partial U_n)
\]
is a weak equivalence. Since 
\[
\beta:KH(G,U;\partial U)\to 
KH(G,U;\partial U_0,\ldots, \partial U_n)
\]
is a weak equivalence, it suffices to show that
\[
K_0(G,U;\partial U)\to KH_0(G,U;\partial U)
\]
is an isomorphism.

We have the commutative diagram, with rows homotopy fiber sequences:
\[
\xymatrix{
K(G,U;\partial U)\ar[r]\ar[d]&K(G,U)\ar[r]\ar[d]&K(G,\partial U)\ar[d]\\
KH(G,U;\partial U)\ar[r]&KH(G,U)\ar[r]&KH(G,\partial U)
}
\]
and 
\begin{align*}
&K_0(G, U)\to KH_0(G, U)\\
&K_1(G, U)\to KH_1(G, U)
\end{align*}
are isomorphisms,  so it suffices to show that the natural maps
\begin{align*}
&K_0(G, \partial U)\to KH_0(G, \partial U)\\
&K_1(G,  \partial U)\to KH_1(G,  \partial U)
\end{align*}
are isomorphisms. This follows from the weakly convergent spectral sequence (\cite[th. 1.3]{WeibelHomo})
\[
E_{p,q}^1=N^pK_q(G,-) \Rightarrow KH_{p+q}(G, -)
\]
and the $K_p(G,-)$ regularity of $\partial U$ for $p\le 1$ (Lemma~\ref{lem:Regularity}).
\end{proof}

\vspace{4mm}

\section{The equivariant homotopy coniveau tower} In this section, we define $G$-equivariant versions of the homotopy coniveau tower defined in \cite{TechLoc}.

Let $G$ be a finite group and $T$ be a $G$-scheme. For a closed $G$-stable subset $W\subset T$
we define
$$G^W(G,T):=\hofib(G(G,T)\rightarrow G(G, T \setminus W)).$$
The localization sequence for equivariant $G$-theory gives us the canonical isomorphism
$$G(G,W)\rightarrow G^W(G, T).$$

We have as well the $K$-theory version:
\[
K^W(G, T):=\hofib(K(G,T)\rightarrow K(G, T \setminus W)).
\]
and the canonical map
\[
K^W(G, T)\to G^W(G, T)
\]
which is a weak equivalence in case $T$ is regular.

We proceed to discuss the functoriality of these constructions. There is always a technical problem occuring at this point, in that the pull-back morphisms on the various exact categories involved are not strictly functorial, but only functorial up to a natural isomorphism satisfying a cocycle condition. Quillen has explained in \cite{AlgKthy1} how to rectify this situation by replacing  the category of locally free coherent sheaves on a scheme $T$ with an equivalent category, for which there exist strictly functorial pull-back maps. We will use this construction throughout, suppressing its explicit mention, so that the $K$-theory and $G$-theory spectra become strictly functorial constructions.

Let $f:T'\to T$ be a $G$-equivariant morphism of  $G$-schemes, $W'\subset T'$ and $W\subset T$ $G$-stable subsets with $f^{-1}(W)\subset W'$. The commutative diagram
\[
\xymatrix{
K(G,T)\ar[r]^{j^*}\ar[d]_{f^*}&K(G,T\setminus W)\ar[d]^{f^*}\\
K(G,T')\ar[r]_{j^*}&K(G,T'\setminus W')}
\]
defines the pull-back map $f^*:K^W(G,T)\to K^{W'}(G,T')$, satisfying the functoriality $(fg)^*=g^*f^*$. The same holds for $G$-theory if $f$ is flat; in fact, if we require that $f$ factors as $p\circ i$, with $p:\P^N\times T \rightarrow T$ the projection and $i:T'\to \P^N\times T$ a $G$-equivariant regular embedding into an open subscheme $U$ of  $\P^N\times T$, then it is not hard to show that $f^*:G^W(G,T)\to G^{W'}(G,T')$ exists in this case as well. Indeed, $p$ is smooth, so it suffices to define $i^*$, and for this, the condition on $i$ and Quillen's resolution theorem show that the inclusion of the $i$-flat $G$-coherent sheafs on 
$\P^N\times T$, $\cM_{G, \P^N\times T,i}$, into $\cM_{G, \P^N\times T}$ induces a weak equivalence
on the $K$-theory spectra. See \cite{TechLoc} for details in the non-equivariant case; the arguments for the equivariant case are exactly the same. A morphism $f$ which admits a factorization $p\circ i$ as above is called an \lci-morphism.

If $T=T'$ and $f=\id$, then $\id^*$ defines the functorial push-forward map
\[
i_{W',W*}:G^W(G,T)\to G^{W'}(G,T)
\]
for $W\subset W'$ $G$-stable closed subsets of $T$.

We recall the cosimplicial $k$-scheme $\Delta^*$, $n\mapsto \Delta^n$;
the cosimplicial structure is  defined by sending an order-preserving map $g:\{0,\ldots, n\}\to \{0,\ldots, m\}$ to the map $g:\Delta^n\to \Delta^m$ with
\[
g^*(t_i):=\sum_{j\in g^{-1}(i)}t_j.
\]

\begin{defn} (1) For $X$ of finite type over $k$, set        
\[
S^{G,X}_{(p)}(r):=\left\{W\subset X\times \Delta^r\vert \begin{array}{c}W \text{ is a closed } G \text{-stable subset} \\
                                                                 \text{ and for all faces } F\subset \Delta^r \text{ we have }\\ 
                                                                 \dim W\cap X\times F \leq p+ \dim F.
                                                             \end{array}\right\}                                                     
\]
(2) For $X$ finite type and locally equi-dimensional over $k$,
\[
S_{G,X}^{(p)}(r):=\left\{W\subset X\times \Delta^r\vert \begin{array}{c}W \text{ is a closed } G \text{-stable subset} \\
                                                                 \text{ and for all faces } F\subset \Delta^r \text{ we have }\\ 
                                                                 \codim_{X\times F}W\cap X\times F\geq p
                                                             \end{array}\right\}
 \]
\end{defn}

The maps $r\mapsto S_{G,X}^{(p)}(r)$,  $r\mapsto S^{G,X}_{(p)}(r)$ define  simplicial sets; if $X$ has pure dimension $d$ over $k$, we have 
\[
S_{G,X}^{(p)}(r)=S^{G,X}_{(d-p)}(r).
\]
We have natural
inclusions $S_{G,X}^{(p+1)}(r)\subset S_{G,X}^{(p)}(r)$,  $S^{G,X}_{(p)}(r)\subset S^{G,X}_{(p+1)}(r)$.

\begin{defn} Let $X$ be a $G$-scheme of finite type over $k$. Let 
\[
G_{(p)}(G,X,r):=\hocolim_{W\in S^{G,X}_{(p)}(r)} G^W(G,X\times\Delta^r).
\]
If $X$ is smooth over $k$, let 
\[
K^{(p)}(G,X,r):=\hocolim_{W\in S_{G,X}^{(p)}(r)} K^W(G,X\times\Delta^r).
\]
\end{defn}

Since $\Delta^n$ is smooth over $k$, the structure morphisms in the cosimplicial scheme $X\times\Delta^*$ are all \lci-morphisms. Thus, the simplicial structure of $S^{G,X}_{(p)}(\_)$ gives us the simplicial spectrum $G_{(p)}(G,X,\_)$. In case $X$ is locally equi-dimensional over $k$ we have as well the simplicial spectrum $G^{(p)}(G,X,\_)$ and if $X$ is smooth over $k$, the simplicial spectrum $K^{(p)}(G,X,\_)$ with term-wise weak equivalence 
\[
K^{(p)}(G,X,\_)\to G^{(p)}(G,X,\_).
\]

We let $\dim X$ denote the maximum of $\dim X_i$ over the irreducible components $X_i$ of $X$. The inclusion $S^{G,X}_{(p)}(r)\subset S^{G,X}_{(p+1)}(r)$  gives rise to the {\em homotopy coniveau tower} of simplicial spectra
\begin{equation}\label{eqn:HCTower}
\dots\rightarrow G_{(p)}(G,X,\_)\rightarrow G_{(p+1)}(G,X,\_)
  \rightarrow \dots \rightarrow G_{(\dim X)}(G,X,\_)
\end{equation}
We denote the layers of this tower by
\[
G_{(p/p-1)}(G,X,\_):= \hocofib (G_{(p-1)}(G,X,\_)\rightarrow G_{(p)}(G,X,\_)).
\]

\begin{rem}   \label{rem:Homotopy}
Since $X\times\Delta^r$ is the final element of $S^{G,X}_{(\dim X)}(r)$,   the canonical map $G(X\times\Delta^r)\to G_{(\dim X)}(X,r)$ is a weak equivalence.  
By homotopy invariance (compare \cite[Cor. 4.2]{ThomasonGroup} and \cite[Th. 5.7]{ThomasonGroup}), the map of the constant simplicial spectrum $G(G,X)$ to  $G_{(\dim X)}(G,X,\_)$ induced by the identity $G(G,X)=G_{(\dim X)}(X,0)$ is a weak equivalence. 
\end{rem}

With that notation we have the following proposition.

\begin{prop}\label{prop:HCSpecSeq}
There is a strongly convergent spectral sequence
\[
E_1^{p,q}=\pi_{-p-q}(G_{(p/p-1)}(G,X,\_))\Rightarrow
G_{-p-q}(G,X).
\]
\end{prop}
\begin{proof}
The spectral sequence is the spectral sequence of the homotopy coniveau tower \eqref{eqn:HCTower} together with the identitfication $G_n(G,X)\cong \pi_n G_{(\dim X)}(G,X,\_)$ given by Remark~\eqref{rem:Homotopy}. Since
\[
G_{(p)}(G,X,r)=0 
\]
for all $r< -p$ and each $G_{(p)}(G,X,r)$ is -1-connected, 
$G_{(p)}(G,X,\_)$ is $-p-1$-connected, and hence the canonical map
\[
G_{(\dim X)}(G,X,\_)\to \hocofib[G_{(p)}(G,X,\_)\to G_{(\dim X)}(G,X,\_)]
\]
is $-p-1$-connected, whence the convergence.
\end{proof}

In case $X$ is smooth over $k$, we have the homotopy coniveau tower for equivariant $K$-theory \begin{equation}\label{eqn:HCKTower}
\dots\rightarrow K^{(p+1)}(G,X,\_)\rightarrow K^{(p)}(G,X,\_)
  \rightarrow \dots \rightarrow K^{(0)}(G,X,\_)
\end{equation}
with layers
\[
K^{(p/p+1)}(G,X,\_):= \hocofib (K^{(p+1)}(G,X,\_)\rightarrow K^{(p)}(G,X,\_)).
\]
and the strongly convergent spectral sequence
\[
E_1^{p,q}=\pi_{-p-q}(K^{(p/p+1)}(G,X,\_))\Rightarrow
K_{-p-q}(G,X).
\]

\section{Local coefficients and the cycle class map} In this section, we define the equivariant cycle complex of Bredon type, $z_q(X,G,*)$, and use the homology of $z_q(X,G,*)$ to define the equivariant higher Chow groups of Bredon type.

\subsection{The local coefficients}
For a $k$-scheme $T$, we let $T_{(p)}$ denote the set of points $t\in T$ whose closure $\overline{\{t\}}$ has $\dim_k\overline{\{t\}}=p$. 

Let $G$ and $X$ be as the the last section.
For any smooth $k$-scheme $Y$ we consider the $G$-scheme $X\times Y$,
where $G$ acts trivially on the second component. Because $G$ is finite
$G$ acts also on the point $(X\times Y)_{(p)}$ of dimension $p$ on $X\times Y$.
For an orbit $[x]\in(X\times Y)_{(p)}/G$ of the $G$-set $(X\times Y)_{(p)}$, take a representative $x\in(X\times Y)_{(p)}$ and let $G_x\subset G$ denote the isotropy subgroup for $x$.   We have the
``local coefficient":
\begin{eqnarray*}
\colim_{U\subset \overline{G\cdot x}} G_0(G,U) & = & G_0(G,\coprod_{y\in G\cdot x} \Spec(\kappa(y))  \\
                                               & = & G_0(G_x,\Spec(\kappa(x))) \\
                                               & = & K_0(G_x,\Spec(\kappa(x))).
\end{eqnarray*}
Here the colimit on the left hand side is taken over all open $G$-invariant subsets $U$.

We have the following functorial properties for the local coefficients.

\begin{prop} \label{prop:coef}
Let $[x]\in (X\times Y)_{(p)}/G$ and let $f:Y'\rightarrow Y$ be a morphism in $\Sm/k$ of pure codimension $q$.  Suppose that 
\[
\dim(1_X\times f)^{-1}(\overline{G\cdot x}) =p-q.
\]
Then there is a well defined pull back morphism
\[
f^*:K_0(G_x,\Spec \kappa(x))\rightarrow \bigoplus_{y\in (X\times Y')_{(p-q)}/G, f(y)\in G\cdot x} K_0(G_y, \Spec \kappa(y))
\]
with$(fg)^*=g^*f^*$ when all three maps are defined.
\end{prop}
\begin{proof} The proof  is a modification of Quillen's proof of Gersten's conjecture \cite{AlgKthy1}.
Let $\cM^1(x)$ be the category of coherent $G$-sheaves on $\overline{G\cdot x}$ whose support contains no generic point, and set $G^{(1)}_n(G, \overline{G\cdot x}):=K_n(\cM^1(x))$. 
Consider the localization sequence
\[
\ldots \to G^{(1)}_0(G,\overline{G\cdot x})\to G_0(G,\overline{G\cdot x})\to
G_0(G, G\cdot x)\to 0.
\]
Since $1\times f$ has finite Tor-dimension, we have the well-defined and functorial map
\[
(1\times f)^*:G_0(G,\overline{G\cdot x})\to 
G_0(G,(1\times f)^{-1}(\overline{G\cdot x}));
\]
it clearly suffices to show that the composition
\begin{multline*}
G^{(1)}_0(G,\overline{G\cdot x})\to G_0(G,\overline{G\cdot x})\\
\xrightarrow{(1\times f)^*}G_0(G,(1\times f)^{-1}(\overline{G\cdot x}))
\to G_0(G,G\cdot y)
\end{multline*}
is the zero map.

If $f:Y'\to Y$ is smooth, and $\cF$ is in $\cM^1(x)$, then $(1\times f)^*\cF$ is in  $\cM^1(y)$, whence the result in this case. Also, if $k$ is finite,  the standard trick of passing to infinite extensions $k\to k_1$, $k\to k_2$, with $k_i$ a union of extensions of degree $\ell_i^r$ for distinct primes $\ell_1, \ell_2$, and using the fact that $p_*\circ p^*=\times m$ on $G_*(G,W)$ for $p:T\to W$ a finite \'etale $G$-map of degree $m$, reduces us to the case of an infinite base-field $k$.

We may replace $Y$ with any open subscheme containing $p_2(G\cdot x)$; as this is a finite set of points of $Y$, we may assume that $Y$ is affine; similarly, we may assume that $Y'$ is affine. Thus we can factor $f$ as a closed embedding $i:Y'\to \A^N\times Y$ followed by the smooth projection $\A^N\times Y\to Y$. This reduces us to the case of a closed embedding of affine schemes in $\Sm/k$. By a similar argument, we may assume that $X$ is affine.

If $i:Y'\to Y$ is such a closed embedding, of codimension say $q$ then we can factor $i$ as a sequence
 of $d$ codimension one closed embeddings
 \[
 Y'=Y_0\to Y_1\to\ldots \to Y_{q-1}\to Y_q=Y
 \]
 with each $Y_i$ smooth in a neighborhood of $p_2(G\cdot y)$. By replacing $Y$ with a suitable neighborhood of   $p_2(G\cdot y)$, we may assume we have a factorization as above with each $Y_i$ smooth over $k$. This reduces us to the case of a codimension one closed embedding of smooth affine $k$-schemes.   
 
 Let $F$ be a $G$-stable closed subset of $\overline{G\cdot x}$ disjoint from $G\cdot x$ and let $\cF$ be a coherent $G$-sheaf supported on $F$. Write $F$ as a union $F=F_0\cup F_1$, where $F_0$ is the union of the irreducible components of $F$ which are disjoint from $G\cdot y$, and $F_1$ is the union of the remaining components.  The localization properties of $G_*(G,-)$  yield a Mayer-Vietoris exact sequence
 \[
 \ldots\to G_0(G,F_0\cap F_1)\to G_0(G,F_0)\oplus G_0(G,F_1)\to G_0(G,F)\to0
 \]
 so we can write $[\cF]=i_{0*}(x_0)+i_{1*}(x_1)$ for elements $x_j\in G_0(G,F_j)$, where $i_j:F_j\to F$ are the inclusions. Since it is clear that $(1\times f)^*i_{0*}(x_0)$ goes to zero in $G_0(G,G\cdot y)$, we may assume that $F=F_1$, i.e. that $F\subset  \overline{G\cdot y}$.

 Suppose $\dim_k Y=n+1$, $\dim_kY'=n$.  Let $S:=p_2(G\cdot y)\subset Y'$. By {\it loc. cit.} there is a morphism $\pi:Y\to \A^n_k$ such that
 \begin{enumerate}
 \item The restriction of $\pi$ to $\bar{\pi}:Y'\to \A^n$ is finite.
 \item $\pi$ is smooth on a neighborhood of $S$ in $Y$.
 \end{enumerate}

 Form the diagram
 \begin{equation}\label{eqn:Quillen}
 \xymatrix{
 Y'\times_{\A^n}Y\ar@<-3pt>[d]_{q_1}\ar[r]^-{q_2}&Y\ar[d]^\pi\\
 Y'\ar[r]_-{\bar{\pi}}\ar@<-3pt>[u]_s&\A^n}
 \end{equation}
 where $s$ is the section to $p_1$ induced by the inclusion $Y'\to Y$. $q_2$ is finite. Since $\pi$ is smooth near $S$, $s(Y')$ is a Cartier divisor in a neighborhood of $q_2^{-1}(S)$. Since $q_2$ is finite, there is a neighborhood $V$ of  $S$ in $Y$ such that $s(Y)$ is principal on $V':=q_2^{-1}(V)$; let $t$ be a defining equation.    
   
 Taking the product of \eqref{eqn:Quillen} with $X$ gives us the diagram of $G$-schemes
  \[
 \xymatrix{
 X\times Y'\times_{\A^n}Y\ar@<-3pt>[d]_{q_1}\ar[r]^-{q_2}&X\times Y\ar[d]^\pi\\
 X\times Y'\ar[r]_-{\bar{\pi}}\ar@<-3pt>[u]_s&X\times \A^n}
 \]
 where we omit the $1_X\times -$ on the morphisms. $T:=p_2^*(t)$ is thus a defining equation for $s(X\times Y')$ over $X\times V'$. 
 
 Now set $U:=X\times V$, $U':=X\times V'$, $D:=s^{-1}(s(X\times Y')\cap U')\subset X\times Y'$, and let $j:U'\to X\times Y'\times_{\A^n}Y$ be the inclusion.  We have the commutative diagram
  \[
 \xymatrix{
 U'\ar[r]^-{q_2'}&U\\
  &D\ar[lu]^s\ar[u]_i}
 \]
 with $s(D)$ defined on $U'$ by $T$. Let $\cG$ be the restriction of $\cF$ to $D$. We have the exact sequence
 \[
0\to q'_{2*}j^* q_1^*\cF\xrightarrow{\times T} q'_{2*}j^*q_1^*\cF\to i_*\cG\to 0;
\]
pulling back to $D$ and noting that $i^*q'_{2*} j^*q_1^*\cF$ is supported in $D\cap\overline{G\cdot y}$ gives the identity in $G_0(G,D\cap\overline{G\cdot y})$
\[
i^*[ i_*\cG]=i^*[q'_{2*} j^*q_1^*\cF]-i^*[q'_{2*} j^*q_1^*\cF]=0.
\]
Restricting to $G\cdot y$ completes the proof.
 \end{proof}

We consider the simplicial set
\[
X_{(p)}^G(r):=\{[x]\in (X\times\Delta^r)_{(p+r)}/G \vert \overline{G\cdot x}\in S_{(p)}^{G,X}(r)\}.
\]
and set
\[
z_p(G,X,r):=\bigoplus_{[x]\in X^G_{(p)}(r)} K_0(G_x,\Spec(\kappa(x))).
\]
 
By Proposition~\ref{prop:coef}, the cosimplicial structure on $r\mapsto\Delta^r$ makes $r\mapsto 
z_p(G,X,r)$ into a simplicial abelian group, denoted $z_p(G,X,\_)$.

\begin{defn} Let $X$ be a finite type $k$-scheme with a $G$-action for a finite group $G$. The
 {\em equivariant cycle complex of Bredon type}, $z_p(X,G,*)$, is the complex associated to the simplicial abelian group $z_p(G,X,\_)$. Define the {\em equivariant higher Chow groups of Bredon type} by
\[
CH_p(G,X,r):=\pi_r(z_p(G,X,\_))=H_r(z_p(X,G,*)).
\]
If $X$ is locally equi-dimensional over $k$, we may  index by codimension, giving us the simplicial abelian group $z^p(G,X,\_)$, the complex $z^p(G,X,*)$ and the codimension $p$ equivariant higher Chow groups $CH^p(X,G,r)$.
\end{defn} 

\subsection{Functorialities} Let  $\rho:L\to F$ be a finite extension of commutative noetherian rings with $G$-action (compatible via $\rho$). $\rho$ induces the exact functor $\rho^*$ (restriction of scalars) from the category of finitely generated $F^{tw}[G]$-modules to finitely generated $L^{tw}[G]$-modules and thereby the map on $G$-theories $\rho^*:G(G, L)\to G(G,F)$. For a $k$-algebra homomorphism $\phi:A\to B$,  we have the natural isomorphism of functors
\[
(\rho\otimes\id)^*\circ(\id\otimes\phi)_*\cong (\id\otimes\phi)_*\circ (\rho\otimes\id)^*.
\]
Thus, if $\phi$ has finite Tor-dimension, the diagram
\[
\xymatrix{
G(G, F\otimes_k A)\ar[r]^-{(\id\otimes\phi)_*}\ar[d]_{(\rho\otimes\id)^*}& 
G(G,F\otimes_kB)\ar[d]^{(\rho\otimes\id)^*}\\
G(G, L\otimes_kA)\ar[r]_-{(\id\otimes\phi)_*}& G(G,L\otimes_kB)
}
\]
commutes.

If $f:Y\to X$ is a proper $G$-equivariant morphism of $k$-schemes with $G$-action, we define the push-forward morphism
\[
f_*(r):z_p(Y,G,r)\to z_p(X,G,r),
\]
by  
\[
[\alpha\in K_0(\kappa(Z),G_z)]\mapsto
[(f\times\id)_*(\alpha)\in K_0(\kappa((f\times\id)(Z)),G_{(f\times\id)(z)})]
\]
 if $Z\to (f\times\id)(Z)$ is generically finite, and sending $\alpha$ to zero if not. By the commutativity of the above diagram, the maps $f_*(r)$ extend to the map of simplicial abelian groups
\[
f_*:z_p(Y,G,\_)\to z_p(X,G,\_),
\]
with $(fg)_*=f_*\circ g_*$ for proper composable morphisms $f,g$.

Similarly, given a flat $G$-equivariant morphism $f:Y\to X$ of relative dimension $d$, we have the pullback map
\[
f^*:z_p(X,G,\_)\to z_{p+d}(Y,G,\_)
\]
with $(fg)^*=g^*f^*$, and we have the compatibility
\[
g^*f_*=f'_*g^{\prime *}
\]
in $G$-equivariant cartesian squares
\[
\xymatrix{
W\ar[r]^{f'}\ar[d]_{g'}&Y\ar[d]^g\\
Z\ar[r]_f&X
}
\]
with $f$ proper and $g$ flat.

\begin{rem} Relying on the localization property of the $z_p(G,X,\_)$, one can   extend the contravariant functoriality of $Y\mapsto z^p(G, X\times_kY, \_)$ from flat morphisms to arbitrary morphisms of smooth $k$-schemes $Y'\to Y$; we will explain this in another paper. However, there seems to be no good pull-back $f^*:z^p(G, X,\_)\to z^p(G, Y, \_)$, compatible with pull-back on equivariant $G$-theory,  for  arbitrary $G$-morphisms $f:Y\to X$, even if both $X$ and $Y$ are smooth over $k$.  We give an example of this phenomenon.

Let $k$ be a field, $Z:=\Spec k[x,y]$, $G:=\ZZ/2$ and let $G$
act on $Z$ via
\[
(1 \text{ mod }2\ZZ)\cdot (x,y):= (-x,-y).
\]
Let $X$ be the line $y=0$, $Y$ the line $x=0$ and $w$
be the intersection of $X$ and $Y$, i. e. $w=\{(0,0)\}.$ Let $y$ be the generic point of $Y$.

We get the following diagram
$$\xymatrix{
K^Y(G,Z) \ar[r] \ar[d]  & K_0(G_y, \kappa(y)) \\
K^w(G,X) \ar[r]              & K_0(G_w, \kappa(w)) }$$
The horizontal arrows are surjective by  localization and the bottom horizontal arrow is an isomorphism. The homotopy property for equivariant $K$-theory implies that the vertical arrow is also an isomorphism, again by the localization theorem.

Since $G_w=G$ acts trivially on $\kappa(w)$, $K_0(G_w, \kappa(w)))$ is isomorphic to the group ring $\ZZ[\ZZ/2]$. However, since $G=G_y$ acts non-trivially on $\kappa(y)$, we have
\[
K_0(G_y, \kappa(y))\cong K_0( \kappa(y)^G)\cong \ZZ.
\]
Thus, there is no arrow on
the right which makes this diagram commutative.

This is related to the fact that the ring struction on $K_0(X,G)$ (given by tensor product over $\cO_X$) does not in general respect the topological filtration, and so one should not expect the equivariant Chow groups $CH^*(G,X,*)$ to have a ring structure. To see an example of this, consider the case $K_0(G,X)$, with $X$ and $G$ as above. It follows from the homotopy property and localization that the image of
\[
K_0^w(G, X)\to K_0(G, X)
\]
is the augmentation ideal of the group ring $\ZZ[\ZZ/2]$. But if $\sigma$ denotes the non-identity element of $\ZZ/2$, we have $(1-\sigma)^2=2(1-\sigma)$. Thus, letting $\ZZ_\sigma$ denote the sign representation of $\ZZ/2$, 
\[
(\ZZ-\ZZ_\sigma)^{\otimes 2}\cong 2(\ZZ-\ZZ_\sigma). 
\]
This shows that $F^1_{top}\cdot F^1_{top}\neq 0$, although $F^2_{top}=0$.

This last point arose during discussions between the first-named author and W. Niziol.
\end{rem}

\subsection{The cycle map}
For each $W\in S_{(p)}^{G,X}(r)$ we have a canonical morphism
\[
\pi_0(G^W(G,X\times\Delta^r))\rightarrow\bigoplus_{[x]\in X^G_{(p)}(r), x\in W} K_0(G_x,\Spec(\kappa(x))).
\]
This defines a morphism
\[
\pi_0(G_{(p)}(G,X,r))\rightarrow\bigoplus_{[x]\in X^G_{(p)}(r)} K_0(G_x,\Spec(\kappa(x))),
\]
which factors through the surjection
\[
\pi_0(G_{(p)}(G,X,r))\rightarrow \pi_0(G_{(p/p-1)}(G,X,r)).
\]
We consider $\pi_0(G_{(p/p-1)}(G,X,n))$ as a  spectrum by using the associated Eilenberg-Maclane spectrum. Noting that  $G_{(p/p-1)}(G,X,n)$ is -1-connected for each $n$, we have the
canonical maps of   spectra 
\[
G_{(p/p-1)}(G,X,n)\to \pi_0(G_{(p/p-1)}(G,X,n)), 
\]
which yield the map of simplicial spectra
\[
G_{(p/p-1)}(G,X,-)\rightarrow \pi_0(G_{(p/p-1)}(G,X,-)).
\]

We similarly consider
$z_p(G,X,\_)$ as a simplicial spectrum. Taking the composition of the maps described above  yields the
{\em cycle map}
\[
cl_p: G_{(p/p-1)}(G,X,\_)\rightarrow z_p(G,X,\_),
\]
defined as a map of simplicial spectra.

\vspace{0.5cm}
Now we can formulate our main result.

\begin{satz}\label{satz:main}
Let $X$ be a   scheme of finite type over a field $k$ with an action of a finite
group $G$. Suppose that $\frac{1}{\sharp G}\in k$. Then the cycle map
\[
cl_p: G_{(p/p-1)}(G,X,\_)\rightarrow z_p(G,X,\_)
\]
is  a weak equivalence for all $p$. 
\end{satz}

This and Proposition~\ref{prop:HCSpecSeq} gives
\begin{cor} There is  a strongly convergent spectral sequence
\[
E_1^{p,q}=CH_p(G,X,-p-q)\Rightarrow G_{-p-q}(G,X).
\]
\end{cor}

For the proof we first reduce in the next section  via localization
techniques to the case where $X$ is a point. In the last section
we discuss the case of a point.

\begin{rem} It is easy to see that the cycle map  and the spectral sequence are natural with respect to proper pushforward and flat pullback.
\end{rem}

\section{Localization and reduction to the point}

The main result of this section is

\begin{satz} \label{satz:localization} Let $G$ be a finite group acting on a finite type $k$-scheme $X$, $i:W\to X$ a $G$-stable closed subscheme with open  complement $j:U\to X$. Then for each $p$, the sequence of simplicial spectra
\[
G_{(p)}(G,W,\_)\xrightarrow{i_*}G_{(p)}(G,X,\_)\xrightarrow{j^*} G_{(p)}(G,U,\_)
\]
is a weak homotopy fiber sequence, and the sequence of complexes
\[
z_p(G,W,*)\xrightarrow{i_*}z_p(G,X,*)\xrightarrow{j^*} z_p(G,U,*)
\]
is isomorphic to a cone sequence in the derived category, i.e., the induced map
\[
(j^*,0):\text{cone}(i_*)\to z_p(G,U,*)
\]
is a quasi-isomorphism.
\end{satz}

\begin{proof} We first consider the $G$-theory sequence. The proof is the same as the proof of the analogous localization theorem in the non-equivariant case \cite[Cor. 8.2]{TechLoc}, replacing the spectra $G_{(p)}(?,\_)$ with the $G$-equivariant versions $G_{(p)}(G,?,\_)$ throughout. For the readers convenience, we give a sketch of the argument.

Let   $S^{G,U^X}_{(p)}(r)$ be the subset of $S^{G,U}_{(p)}(r)$ consisting of those $W\subset U\times\Delta^r$ with closure $\overline{W}\subset X\times\Delta^r$ in $S^{G,U}_{(p)}(r)$. These form a simplicial subset of $r\mapsto S^{G,U}_{(p)}(r)$; let 
\[
G_{(p)}(G,U^X,r):=\hocolim_{W\in S^{G,U^X}_{(p)}(r)} G^W(G,X\times\Delta^r),
\]
forming the simplicial spectrum $G_{(p)}(G,U^X,\_)$. The inclusions
$S^{G,U^X}_{(p)}(r)\subset S^{G,U}_{(p)}(r)$ induce the map of simplicial spectra
\[
\iota:G_{(p)}(G,U^X,\_)\to G_{(p)}(G,U,\_),
\]
the localization sequence of $G$-equivariant $G$-theory gives the weak homotopy fiber sequence
\[
G_{(p)}(G,W,r)\xrightarrow{i_*}G_{(p)}(G,X,r)\xrightarrow{j^*} G_{(p)}(G,U^X,r)
\]
for each $r$, and thus the weak homotopy fiber sequence
\[
G_{(p)}(G,W,\_)\xrightarrow{i_*}G_{(p)}(G,X,\_)\xrightarrow{j^*} G_{(p)}(G,U^X,\_).
\]
Thus, we must show that $\iota$ is a weak equivalence. Since both simplicial spectra are -1-connected, the Hurewicz theorem tells us that it suffices to show that $\iota$ is a homology isomorphism. 

For a spectrum $E$, let $\ZZ E$ be a functorial model for a chain complex with homology groups the homology of $E$. In particular, for $W\subset U\times \Delta^r$, we have the complex 
$\ZZ G^W(G,U\times\Delta^r)$
representing the homology of $G^W(G,U\times\Delta^r)$. Taking the limit of $W\in
S^{G,U}_{(p)}(r)$ or in $S^{G,U^X}_{(p)}(r)$ gives us the complexes $\ZZ G_{(p)}(G,U,r)$ and 
$\ZZ G_{(p)}(G,U^X,r)$ computing the homology of $G_{(p)}(G,U,r)$ and $G_{(p)}(G,U^X,r)$.

For $W\subset U\times\Delta^r$, let $W_n\subset U\times\Delta^n$ be the union of $(\id\times g)^{-1}(W)$, as $g:\Delta^n\to\Delta^r$ runs over structure morphisms for the cosimplicial scheme $\Delta$. Using the usual alternating sum of the  pullback by coface maps $\id\times\delta^r_i:U\times\Delta^r\to U\times\Delta^{r+1}$, we form the double complex $n\mapsto
\ZZ G^{W_n}(G,U\times\Delta^n)$ and denote the associated total complex by $\ZZ G^W(\ZZ U\times\Delta^*)$. Thus the limit of the complexes $\ZZ G^W(\ZZ U\times\Delta^*)$ over $W\in 
S^{G,U}_{(p)}(r)$ or in $S^{G,U^X}_{(p)}(r)$, $r=1, 2,\ldots$,  computes the homology of $G_{(p)}(G,U,\_)$ and $G_{(p)}(G,U^X,\_)$. We denote the  limits of these complexes by 
$\ZZ G_{(p)}(G,U)^*$ and $\ZZ G_{(p)}(G,U^X)^*$, respectively.  It thus suffices to show that
\[
\iota_\ZZ: \ZZ G_{(p)}(G,U^X)^*\to  \ZZ G_{(p)}(G,U)^*
\]
is a quasi-isomorphism.

For 
$W\in S^{G,U}_{(p)}(r)$, $W'\in S^{G,U^X}_{(p)}(r)$, let
\[
\iota_W:\ZZ G^W(G,U\times\Delta^*)\to \ZZ G_{(p)}(G,U)^*
\]
 and 
 \[
\iota_{W'}^X:\ZZ G^W(G,U\times\Delta^*)\to \ZZ G_{(p)}(G,U^X)^*
\]
be the canonical maps.

Next, we construct another pair of complexes which approximate $\ZZ G_{(p)}(G,U)^*$ and $\ZZ G_{(p)}(G,U^X)^*$. For this, fix an integer $N\ge0$. Let $\partial\Delta^N_i\subset \Delta^N$ be the subscheme defined by $t_i=0$; for $I\subset \{0,\ldots, N\}$ let $\partial\Delta^N_I$ be the face $\cap_{i\in I}\partial\Delta^N_i$. For $I\supset J$, let $i_{J,I}:\Delta^N_I\to\Delta^N_J$ be the inclusion.

Let $\ZZ\Sm/k$ be the additive category generated by $\Sm/k$, i.e., for connected $X$, $Y$, $\Hom_{\ZZ\Sm/k}(X,Y)$ is the free abelian group on the set of morphisms  $\Hom_{\Sm/k}(X,Y)$, and disjoint union becomes direct sum. We will construct objects in the category of complexes $C(\ZZ\Sm/k)$.

Form the complex $(\Delta^N,\partial\Delta^N)^*$ be the complex which is $\oplus_{I,\ |I|=n}\partial\Delta^N_I$ in degree $-n$, and with differential 
\[
d^{-n}:(\Delta^N,\partial\Delta^N)^{-n}\to (\Delta^N,\partial\Delta^N)^{-n+1}
\]
given by $d^{-n}:=\prod_{I,\ |I|=n}d^{-n}_I$, where
\[
d^{-n}_I:\partial\Delta^N_I\to \oplus_{J,\ |J|=n-1}\partial\Delta^N_J
\]
is the sum
\[
d^{-n}_I:=\sum_{j=1}^n i_{I\setminus\{i_j\},I},
\]
where $I=(i_1,\ldots, i_n)$, $i_1<\ldots<i_n$.

We also have the complex $\ZZ\Delta^*$, which is $\Delta^n$ in degree $n$, with differential the usual alternating sum of coboundary maps.

The identity map on $\Delta^N$ extends to a map of complexes
\[
\Phi^N:\ZZ\Delta^*\to (\Delta^N,\partial\Delta^N)[-N];
\]
the maps in degree $r<N$ are all $\pm\id_{\Delta^r}$. We can take the product of   this construction with $U$, giving us the complex 
$U\times(\Delta^N,\partial\Delta^N)$ and the map  of complexes
\[
 \Phi^N:U\times\ZZ\Delta^*\to U\times(\Delta^N,\partial\Delta^N)[-N]
\]

 For $W\in S^{G,U}_{(p)}(N)$, form the complex $\ZZ G^W(G,U\times(\Delta^N,\partial\Delta^N))$ by taking $\oplus_{I, |I|=n} \ZZ G^{W_{N-n}}(G,U\times\Delta^{N-n})$ in degree $-n$, using the differentials in $U\times(\Delta^N,\partial\Delta^N)$ to form a double complex and then taking the total complex. We thus have the map of complexes
\[
\Phi^{N*}_W: \ZZ G^W(G,U\times(\Delta^N,\partial\Delta^N))[-N]\to 
\ZZ G^W(G,U\times\Delta^*).
\]
One shows that $\Phi^{N*}_W$ induces a homology isomorphism in degrees $<N$ (see
\cite[Lemma 2.6]{TechLoc}).

Take $W\in S^{G,U}_{(p)}(N)$. The main result of  \cite{TechLoc}, Theorem 1.9,  gives a map of complexes
\[
\Psi_W:U\times\ZZ\Delta^*\to U\times(\Delta^N,\partial\Delta^N)[-N]
\]
and a degree -1 map
\[
H_W:U\times\ZZ\Delta^*\to U\times(\Delta^N,\partial\Delta^N)[-N]
\]
with the following properties:
\begin{enumerate}
\item $dH_W=\Psi_W-\Phi^N$.
\item Write $\Psi_W$ as a sum 
\[
\Phi_W=\sum_{i=0}^N\sum_{I,j\ |I|=i} n_{ij}\psi_{ijI}
\]
with $\psi_{ijI}:\Delta^{N-i}\to \partial\Delta^N_I=\Delta^{N-i}$ maps in $\Sm/k$. Then $\psi_{ijI}^{-1}(W_{N-i})$ is in $S^{G,U^X}_{(p)}(N-i)$.
\item Write $H_W$ as a sum 
\[
H_W=\sum_{i=0}^N\sum_{I,j\ |I|=i} n_{ij}H_{ijI}
\]
with $H_{ijI}:\Delta^{N-i+1}\to \partial\Delta^N_I=\Delta^{N-i}$ maps in $\Sm/k$. Then $H_{ijI}^{-1}(W_{N-i})$ is in $S^{G,U}_{(p)}(N-i+1)$. If $W'\subset W_{N-i}$ is in $S^{G,U^X}_{(p)}(N-i)$, then 
$H_{ijI}^{-1}(W')$ is in $S^{G,U^X}_{(p)}(N-i+1)$.
\end{enumerate}

Thus $\Psi_W$ induces the map of complexes
\[
\Psi_W^*: \ZZ G^W(G,U\times(\Delta^N,\partial\Delta^N))[-N]\to
 \ZZ G_{(p)}(G,U^X)^*
 \]
 and $H_W$ gives a degree 1 map
 \[
 H_W^*: \ZZ G^W(G,U\times(\Delta^N,\partial\Delta^N))[-N]\to
 \ZZ G_{(p)}(G,U)^*
 \]
 with 
 \[
 dH_W^*=\iota_\ZZ\circ\Psi_W^*-\iota_W\circ \Phi^{N*}_W
 \]
 Furthermore, if $W'\subset W$ is in $S^{G,U^X}_{(p)}(N)$, then $H_W$ gives a
  degree 1 map
 \[
H_W^{X*}: \ZZ G^{W'}(G,U\times(\Delta^N,\partial\Delta^N))[-N]\to
 \ZZ G_{(p)}(G,U^X)^*
 \]
 with 
 \[
 dH_W^{X*} = \Psi_W^*-\iota_{W'}^X\circ \Phi^{N*}_{W'}.
 \]
 Since $\Phi^{N*}_W$ is a homology isomorphism in degrees $<N$ and $\ZZ G_{(p)}(G,U)^*$ and
 $\ZZ G_{(p)}(G,U^X)^*$ are the limits of $\ZZ G^W(G,U\times\Delta^*)$ and $\ZZ G^{W'}(G,U\times\Delta^*)$, respectively, this shows that $\iota_\ZZ$ is a quasi-isomorphism, completing the proof.
 
 For the sequence of cycle complexes, let $z_p(G,U^X,r)$ be the subgroup of $z_p(G,U,r)$ generated by the irreducible codimension $p$ closed subsets $W\subset U\times\Delta^r$ with $W\in S^{G,U^X}_{(p)}(r)$. This forms the subcomplex $z_p(G,U^X,*)$ of $z_p(G,U,*)$ and gives us 
the term-wise exact sequence of complexes
 \[
0\to z_p(G,W,*)\xrightarrow{i_*}z_p(G,X,*)\xrightarrow{j^*} z_p(G,U^X,*)\to0.
\]
Thus, we need to show that the inclusion
\[
z_p(G,U^X,*)\to z_p(G,U,*)
\]
is a quasi-isomorphism. The proof is now exactly the same as the case of $G$-theory, except that we can avoid the use of the Hurewicz theorem by working directly with the complexes $z_p(G,?,*)$ instead of passing to complexes representing the homology of the simplicial abelian group $z_p(G,?,\_)$.
 \end{proof}
 
 \begin{cor} \label{cor:localization} With the hypotheses and notations as in Theorem~\ref{satz:localization}, the sequence
\[
G_{(p/p-1)}(G,W,\_)\xrightarrow{i_*}G_{(p/p-1)}(G,X,\_)\xrightarrow{j^*} G_{(p/p-1)}(G,U,\_)
\]
 is a weak homotopy fiber sequence for all $p$.
 \end{cor}
 
 \begin{proof} This follows directly from Theorem~\ref{satz:localization}, the naturality of $i_*$ and $j^*$ with respect to change of $p$, and the Quetzalcoatl lemma.
 \end{proof}
 
 \begin{cor} \label{cor:FiniteReduction} Suppose that Theorem~\ref{satz:main} is true for all fields $k$ with $\frac{1}{\sharp G}\in k$ and with $X=\Spec K$, where $K$ is a finite  extension of $k$ with a $G$-action such that $G$ acts trivially on $k$. Then Theorem~\ref{satz:main} is true for all all fields $k$ with $\frac{1}{\sharp G}\in k$.
 \end{cor}
 
 \begin{proof}We prove Theorem~\ref{satz:main} by  induction on $\dim_kX$, the case of dimension 0 being true by hypothesis.

 We have already remarked that the cycle map
  \[
cl_p: G_{(p/p-1)}(G,X,\_)\rightarrow z_p(G,X,\_)
\]
is natural with respect to proper pushforward   and pullback with respect to flat maps. Thus, for each $G$-stable closed subset $W\subset X$, we have the commutative diagram
\[
\xymatrix{
G_{(p/p-1)}(G,W,\_)\ar[r]^{i_*}\ar[d]_{cl_p^W}&G_{(p/p-1)}(G,X,\_)\ar[r]^{j^*}\ar[d]_{cl_p^X}& G_{(p/p-1)}(G,U,\_)\ar[d]_{cl_p^U}\\
z_p(G,W,\_)\ar[r]_{i_*}&z_p(G,X,\_)\ar[r]_{j^*} &z_p(G,U,\_)
}
\]
By Theorem~\ref{satz:localization} and Corollary~\ref{satz:localization} the rows are weak homotopy fiber sequences; by our induction hypothesis, $cl_p^W$ is a weak equivalence, so $cl_p^X$ is a weak equivalence if and only if $cl_p^U$ is.
Taking the limit over all open dense $G$-stable $U\subset X$ reduces us to showing that 
\[
cl_p^{\kappa(X)}:G_{(p)}(G,\kappa(X),\_)\to z_p(G,\kappa(X),\_)
\]
is a weak equivalence. Breaking up $\kappa(X)$ into a product of fields reduces us to the case of irreducible $X$.

Note that $G_{(p)}(G,\kappa(X),\_)$, $z_p(G,\kappa(X),\_)$ and $cl_p$ do not depend on the choice of ``constants" $k\subset \kappa(X)$, so we may replace $k$ with the invariant subfield $\kappa(X)^G$. Since the extension $\kappa(X)^G\subset \kappa(X)$ is finite, $ cl_p^{\kappa(X)}$ is a weak equivalence by hypothesis, completing the proof.
\end{proof}

\section{The case of the point}
Now we consider the case of a point. Let $X=\Spec(K)$, $K\supset k$ a finite   field extension of $k$ with a $G$-action,  with $\frac{1}{\sharp G}\in K$ such that $G$ acts trivially on $k$.
Let $\hat{k}:=K^G$ be the subfield of $K$ which is fixed under the operation of $G$. Then the field extension
$K/\hat{k}$ is finite and Galois; as in the proof of Corollary~\ref{cor:FiniteReduction}, we may replace $k$ with $\hat{k}$. Changing notation, we assume that $k$ is the fixed subfield of $K$ under $G$ . We consider the following functor:
$$
\begin{array}{cccc}
  E: & \Sm/k^{op} & \rightarrow & \Spt \\
     & Y    & \mapsto     & K(G,Y\otimes_k K)
\end{array}
$$ 
This clearly defines a presheaf of spectra on $\Sm/k$.

We want to show that this presheaf satisfies the axioms
of \cite{LeHo}. We first recall some notations from \cite{LeHo}.

Let $v(n):=\{v_0(n),\ldots, v_n(n)\}$ be the vertices of $\Delta^n$, where $v_i(n)$ is the point $t_j=0$, $j\neq i$, $t_i=1$.  For a field $F$, let $\cO_{\Delta^n_F,v(n)}$ be the semi-local ring of $v(n)$ in $\Delta^n_F$ and set $\Delta^n_{0,F}:=\Spec\cO_{\Delta^n_F,v(n)}$. The $\Delta^n_{0,F}$ form a cosimplicial subscheme of $n\mapsto \Delta^n_F$. We extend this notation to $F$ a product of fields in the evident manner.

We let $\partial_i\Delta^n_{0,F}:=\partial_i\Delta^n_{F}\cap \Delta^n_{0,F}$, i.e. $\partial_i\Delta^n_{0,F}$ is the closed subscheme of $\partial\Delta^n_{0,F}$ defined by $t_i=0$, and set  $\partial\Delta^n_{0,F}:=\cup_{i=0}^n\partial_i\Delta^n_{0,F}$. We let $\partial_*\Delta^n_{0,F}$ denote the set of components of
$\partial\Delta^n_{0,F}$, 
\[
\partial_*\Delta^n_{0,F}:=\{\partial_0\Delta^n_{0,F},\ldots,\partial_n\Delta^n_{0,F}\}
\]
and denote, e.g., the relative $K$-theory $K(\Delta^n_{0,F};\partial_0\Delta^n_{0,F},\ldots,\partial_n\Delta^n_{0,F})$ by
$K(\Delta^n_{0,F};\partial_*\Delta^n_{0,F})$.

Let   $E:\Sm/k^\op\to \Spt$ be presheaf of spectra on $\Sm/k$. We define the presheaf of spectra $\Omega_TE:\Sm/k^\op\to \Spt$ by
\[
\Omega_TE(X):=E^{X\times0}(X\times\A^1).
\]

From \cite[Section 6]{LeHo}, we have the notion of a {\em well-connected} functor
$E:\Sm/k^\op\to \Spt$. From the definition  \cite[Definition 6.1.1]{LeHo} and \cite[Proposition 6.3.3]{LeHo}, $E$ is well-connected if $E$ satisfies:

\begin{enumerate}
\item Homotopy invariance:\\
For each $X\in \Sm/k$ the
map $p^*:E(X)\rightarrow E(\mathbb{A}^1\times X)$ is
an weak equivalence. 
\item Nisnevich excision:\\
Let $f: X'\rightarrow X$ be an \'etale morphism in $Sm/k$,
and let $W\subset X$ be a closed subset. Suppose that
$f$ restricts to an isomorphism $f^{-1}(W)\rightarrow W$.
Then  
\[
f^*:E^W(X)\rightarrow E^{f^{-1}(W)}(X')
\]
is a weak equivalence.
\item Finite descent:\\
Suppose the base-field $k$ is finite. Let $k'\supset k$ be a finite Galois extension of $k$ with group $H$. Consider the presheaf $E_{k'}$ on $\Sm/k$ defined by $E_{k'}(Y):=E(Y\times_kk')$. $H$ acts on $E_{k'}$; let $E_{k'}^H(Y)$ be a functorial model for the $H$-homotopy fixed point spectrum, giving the presheaf $E_{k'}^H$ and the natural transformation $\theta_{k'}:E\to E_{k'}^H$. Then $\theta_{k'}$ is a weak equivalence after inverting $\sharp H$.

\item Well-connectedness (1): \\
For a closed $W$ in a smooth $X$ the
spectrum $E^W(X)$ is $-1$ connected.  
\item Well connectedness (2): \\
For all $n\ge1$, $d\ge0$, and all finite generated fields $F$ over $k$,
\[
\pi_0[(\Omega_T^d E)(\Delta_{0,F}^n, \partial_*\Delta_{0,F}^n)]=0.
\]
\end{enumerate}

\begin{prop} \label{prop:main} Let $K$ be a field with $G$-action, $k:=K^G$. Suppose that $\frac{1}{\sharp G}\in K$. Suppose further that the presheaf $E: \Sm/k^{op} \to \Spt $,
\[
E(X) :=K(G,X\otimes_k K),
\]
is well-connected. Then for  all $X\in \Sm/k$,  the cycle map
\[
cl_p: G_{(p/p-1)}(G,X\otimes_kK,\_)\rightarrow z_p(G,X\otimes_kK,\_)
\]
is  a weak equivalence for all $p$. 
\end{prop}

\begin{proof} We use the notation of \cite[Section 5]{LeHo}, except that we index with respect to dimension rather than codimension to maintain the conventions used here.

Since $K/k$ is finite galois with galois group $G$ the obvious map
$$ S^X_{(p)}(r) \rightarrow S^{G,X\otimes_k K}_{(p)}(r)$$
is a bijection and therefore we can identify
\begin{align*}
E_{(p)}(X,-) &= G_{(p)}(G, X\otimes_kK,-) \\
E_{(p/p+1)}(X,-) &= G_{(p/p+1)}(G, X\otimes_kK,-) \\
z_p(X;E,-)       &= z_p(G,X\times_kK,-). 
\end{align*}
Further the maps $cl_p$ are compatible with these identifications.

Since $E$ is well-connected, it follows from 
\cite[Corollary 5.3.2]{LeHo} that for $X\in\Sm/k$, the simplicial spectrum 
\[
E_{(p/p-1)}(X,\_)=G_{(p/p-1)}(G,X\otimes_kK,\_)
\]
 is weakly equivalent to  the simplicial spectrum 
$E^{\text{s.l.}}_{(p/p-1)}(X,\_)$. In addition, there is for each $n$ a weak equivalence
\[
cl^{\text{s.l.}}_{p,n}:E^{\text{s.l.}}_{(p/p-1)}(X,n)\to z_p(G,X,n)
\]
(more precisely, to the Eilenberg-Maclane spectrum associated to the abelian group $z_p(G,X,n)$). 

The argument of \cite[Theorem 6.4.1]{LeHo}, repeated word for word, shows that $cl^{\text{s.l.}}_{p,n}$ induces a weak equivalence of simplicial spectra
\[
cl^{\text{s.l.}}_p:E^{\text{s.l.}}_{(p/p-1)}(X,n)\to z_p(G,X,n)
\]
and that composing $cl^{\text{s.l.}}_p$ with the weak equivalence $G_{(p/p-1)}(G,X\otimes_kK,\_)\to
E^{\text{s.l.}}_{(p/p-1)}(X,\_)$ yields the map $cl_p$ (in the stable homotopy category).
\end{proof}

Thus, with the help of Corollary~\ref{cor:FiniteReduction}, to finish the proof of Theorem~\ref{satz:main} we need only show:

\begin{prop} The presheaf $E: \Sm/k^{op} \to \Spt $,
\[
E(Y) =K(G,Y\otimes_k K),
\]
is well-connected.
\end{prop}

\begin{proof} 

We need to show that $E$ satisfies:

\begin{enumerate}
\item[(i)] Homotopy invariance. \\
This follows immediately from
the homotopy invariance for equivariant $K$-theory.
\item[(ii)] Nisnevich excision. \\
This follows from
the fact that the equivariant localization theorem
implies $E^W(X)=G(G,W\otimes_k K)$. 
\item[(iii)] Finite descent. \\
Let $K'=K\otimes_kk'$ and let $\pi:\Spec K'\to \Spec K$ be the projection. Finite descent follows from the existence of the natural pushforward map
 \[
 \pi_*:K(G,Y\otimes_kK')\to K(G,Y\otimes_kK)
 \]
 with $\pi_*\circ\pi^*=\times\sharp H$ on the $K$-groups  $K_n(G,Y\otimes_kK)$.
\item[(iv)] Well-connectedness (1). \\
Since the equivariant $G$-theory spectrum is -1-connected,
 $E^W(X)$ is $-1$ connected  by the
localization property, as in (ii).
\item[(v)] Well connectedness (2). \\
\end{enumerate}

To prove (v), we first note that, since $k\to K$ is finite and separable, $F\otimes_kK$ is a finite product of fields and we have
\begin{align*}
&\Delta_{0,F}^n\times_kK\cong \Delta_{0,F\otimes_kK}^n\\
&\partial_i\Delta_{0,F}^n\times_kK\cong \partial_i\Delta_{0,F\otimes_kK}^n
\end{align*}
Also, the projective bundle formula for $K(G,\;_)$ \cite{ThomasonGroup} shows that $\Omega_TK(G,\;_)=
K(G,\;_)$. Thus (v) follows from

\begin{lem}
Let $F$ be a finitely generated field extension of $k$. Then 
$$K_0(G,\Delta_{0,F\otimes_kK}^n;\partial_*\Delta_{0,F\otimes_kK}^n)=0.$$
\end{lem}
\begin{proof} Let $A=F\otimes_kK$. $G$ acts transitively on the irreducible components of $\Spec A$. Fix one component $x:=\Spec \kappa(x)$, and let $G_x$ be the isotropy group of $x$. Then
\[
K_0(G,\Delta_{0,F\otimes_kK}^n;\partial_*\Delta_{0,F\otimes_kK}^n)=
K_0(G_x,\Delta_{0,\kappa(x)}^n;\partial_*\Delta_{0,\kappa(x)}^n).
\]
In addition, letting $k_x\subset K$ be the fixed field of $G_x$, we have
\[
F\otimes_{k_x}K=\kappa(x).
\]
Thus, changing notation, we may assume that $F\otimes_kK$ is a field $L$ and $F$ is the $G$-fixed subfield of $L$.

We note that $\Delta^n_{0,L}$ is an intersection of $G$-stable affine open subschemes $U$ of $\Delta^n_L$. Thus by Proposition~\ref{homfibseq}, 
the map
\[
\alpha:K_0(G,\Delta_{0,L}^n;\partial\Delta_{0,L}^n)\to
K_0(G,\Delta_{0,L}^n;\partial_*\Delta_{0,L}^n)
\]
is an isomorphism.

We have the exact sequence
\begin{multline*}
K_1(G,\Delta_{0,L}^n)\rightarrow K_1(G,\partial\Delta_{0,L}^n)
\rightarrow K_0(G,\Delta_{0,L}^n;\partial\Delta_{0,L}^n)\\ \rightarrow 
K_0(G,\Delta_{0,L}^n) \rightarrow K_0(G, \partial\Delta_{0,L}^n)
\end{multline*}

Let $R:=\cO_{\Delta^n_{0,F},v(n)}$, $\bar{R}:=R/(\prod_{i=0}^nt_i)$. Since $F\to L$ is finite, we have
\begin{align*}
 &\Delta_{0,L}^n=\Spec R\otimes_FL\\
 &\partial\Delta_{0,L}^n=\Spec \bar{R}\otimes_FL
 \end{align*}
 Also, by Lemma~\ref{twiequ}
 \begin{align*}
 K(G,\Delta_{0,L}^n)=K(R\otimes_F L^{tw}[G])\\
  K(G,\partial\Delta_{0,L}^n)=K(\bar{R}\otimes_F L^{tw}[G])
  \end{align*}
By the lemma below, $R\otimes_F L^{tw}[G]$ and $\bar{R}\otimes_F L^{tw}[G]$ are semi-local rings, and the surjection $R\otimes_F L^{tw}[G]\to \bar{R}\otimes_F L^{tw}[G]$ induces a bijection on the (finite) sets of maximal two-sided ideals.  It follows easily from this that $K_1(R\otimes_F L^{tw}[G])\to K_1(\bar{R}\otimes_F L^{tw}[G])$ is surjective. 

Using the notations of the lemma below, we have
\begin{align*}
R\otimes_FL^{tw}[G]=\prod_{i=1}^r M_{n_i}(R\otimes_FD_i)\\
\bar{R}\otimes_FL^{tw}[G]=\prod_{i=1}^r M_{n_i}(\bar{R}\otimes_FD_i).
\end{align*}
In addition, there is a finite separable field extension $F\subset F_i$ such that $R\otimes_FD_i$ is an Azumaya algebra over $R\otimes_FF_i$ and  $\bar{R}\otimes_FD_i$ is an Azumaya algebra over $\bar{R}\otimes_FF_i$. Since $R\otimes_FF_i$ is integral and semi-local, this implies that each projective module over $R\otimes_FD_i$ is free; by Morita equivalence we have
\[
K_0(M_{n_i}(R\otimes_FD_i))=\ZZ,
\]
generated by the class of $(R\otimes_FD_i)^{n_i}$. This easily implies that the map
\[
K_0(M_{n_i}(R\otimes_FD_i))\to K_0(M_{n_i}(\bar{R}\otimes_FD_i))
\]
is injective, hence $K_0(R\otimes_F L^{tw}[G])\to K_0(\bar{R}\otimes_F L^{tw}[G])$ is injective.

Thus $ K_0(G,\Delta_{0,L}^n;\partial\Delta_{0,L}^n)=0$, completing the proof.
\end{proof}

\begin{lem} Let $L$ be a field with an action of a finite group $G$, $F\subset L$ the fixed subfield of $L$. Suppose that $\frac{1}{\sharp G}\in F$.
Let  $R$ be a noetherian commutative reduced semi-local $F$-algebra
with the property that for each maximal ideal $m\subset R$ we have
$R/m\simeq F$. Then 
\[
R\otimes_FL^{tw}[G]=\prod_{i=1}^r M_{n_i}(R\otimes_FD_i)
\]
where each $D_i$ is a central division algebra over $F_i$ for some finite separable field extensions $F\subset F_i$, and $M_{n_i}$ denotes the $n_i$ by $n_i$ matrix algebra. In addition, the $D_i$, $F_i$, $n_i$ and the integer $r$ are independent of $R$, and the center of $M_{n_i}(R\otimes_FD_i)$ is $R\otimes_FF_i$.
\end{lem}

\begin{proof}
Since $\frac{1}{\sharp G}\in F$,  $L^{tw}[G]$ is a semi-simple, separable algebra over $F$ containing $F$ in its center. Thus, we have $L^{tw}[G]=\prod_{i=1}^{n}
M_{n_i}(D_i)$ where $D_i$ are division rings containing $F$ in their centers, $n_i\in\NN$.  Let $F_i$ be the center of $D_i$. Since $F\to L^{tw}[G]$ is a finite separable extension, each $F_i$ is finite and separable over $F$.

 We claim that $R\otimes_FF_i$ is the center of $R\otimes_FD_i$. Indeed, the center $\cZ_i$ of 
$R\otimes_FD_i$ is a finitely generated projective $R$-module containing $R\otimes_FF_i$. By our assumption on the maximal ideals of $R$, we have $R/m\otimes_R\cZ_i=F_i=R/m\otimes_RR\otimes_FF_i$ for all maximal ideals $m\subset R$, so by Nakayama's lemma, 
$R\otimes_FF_i=\cZ_i$.
\end{proof}

This completes the proof of the proposition, and the proof of Theorem~\ref{satz:main}
\end{proof}

\bibliographystyle{alpha}
\bibliography{Literatur.bib}

\end{document}